\documentclass[a4paper,10pt, reqno, english]{amsart}

\usepackage[utf8]{inputenc}
\usepackage{url,paralist}
\usepackage{caption}
\usepackage{subcaption}
\usepackage{makecell} 
\usepackage[colorlinks=true,urlcolor=blue,citecolor=magenta]{hyperref}
\usepackage{color}
\usepackage{float}
\usepackage{dsfont}
\usepackage{tikz}
\usepackage{amsfonts,amssymb,enumerate}
\usepackage{graphicx}
\usepackage{amsmath,bm,amsthm}
\usepackage{mathtools}
\usetikzlibrary{arrows.meta}
\usepackage[arrow,curve,matrix,tips,2cell]{xy}  
  \SelectTips{eu}{10} \UseTips
  \UseAllTwocells

\usepackage[export]{adjustbox}
\usepackage{capt-of}

\theoremstyle{plain}
\newtheorem{theorem}{Theorem}[section]
	
\newtheorem{lemma}[theorem]{Lemma}
\newtheorem{proposition}[theorem]{Proposition}
\newtheorem{corollary}[theorem]{Corollary}

\newtheorem{problem}[theorem]{Problem}

\theoremstyle{definition}
\newtheorem{definition}[theorem]{Definition}

 \usepackage{setspace}

{\catcode`@=11\global\let\c@equation=\c@theorem}

\newcommand{\IR}{{\mathbb R}} 
\newcommand{\IF}{{\mathbb F}} 

\newcommand{\id}{\operatorname{id}} 
\newcommand{\Gr}[2]{ G_{#1}( \mathbb{R}^{#2}) } 

\newcommand{\depth}{\mathrm{depth}}
\newcommand{\Fl}{\mathrm{Fl}}

\newcommand{\xySquare}[9] {	 
	\xymatrix@C=#9mm{
		#1 \ar[r]^{#2} \ar[d]_{#4} & #3 \ar[d]^{#5}  \\
		#6\ar[r]^{#7} & #8 }
}

\newcommand{\FiberBundle}[4] {	 
	\xymatrix{
		#1 \; \ar@{^{(}->}[r] & #2 \ar[d]^{#3}  \\
		& #4 }
}

\newcommand{\version}[1]{ 
	\begin{center} 
		last edited on #1\\
		last compiled on \today\\
		name of texfile: \jobname
	\end{center}
}

\title[ A Center Transversal Theorem for mass assignments]
{ A Center Transversal Theorem for mass assignments }

\author{Omar Antol\'in Camarena}
\address{
	\hfill\break Instituto de Matem\'aticas, Universidad Nacional Aut\'onoma de M\'exico,  C.U.  \\
	\hfill\break CDMX 04510, Mexico.
}
\email{omar@matem.unam.mx}

\author{Jaime Calles Loperena}
\address{
	\hfill\break Instituto de Matem\'aticas, Universidad Nacional Aut\'onoma de M\'exico,  C.U.  \\
	\hfill\break CDMX 04510, Mexico.
}
\email{calles@im.unam.mx}

\begin{document}

\date{\today}
\maketitle

\begin{abstract}

In this paper, based on the ideas of Blagojevi\'c, Karasev \& Magazinov, we consider an extension of the center transversal theorem to mass assignments with an improved Rado depth. In particular we substitute the marginal of a measure by a more general concept called a mass assignment over a flag manifold. Our results also allow us to solve the main problem proposed by Blagojevi\'c, Karasev \& Magazinov in a linear subspace of lower dimension, as long as it is contained in a high-dimensional enough ambient space.

\end{abstract}

\section{ \Large Introduction}
\label{sec: Intro}

Motivated by the previous work of Blagojevi\'c, Karasev \& Magazinov \cite{blagojevic2018center}, Schnider \cite{schnider2020ham}, and Axelrod-Freed \& Sober\'on \cite{axelrod2022bisections}, we consider an extension of the classical center transversal theorem to mass assignments. For that purpose, we will start by establishing the terminology that we will use.

\medskip
Let $v\in S^{d-1}$ be a unit vector in $\mathbb{R}^d$, and let $a \in \mathbb{R}$. An {\em oriented affine hyperplane}  $H_{v,a}:=\{x \in \mathbb{R}^d \colon \langle x,v\rangle =a\}$ in $\IR^d$ determines two closed half-spaces denoted by  
\[
	H_{v,a}^{0}:=\{x \in \mathbb{R}^d \colon\langle x,v\rangle \geq a\}\qquad\text{and}\qquad
		H_{v,a}^{1}:=\{x \in \mathbb{R}^d \colon \langle x,v\rangle \leq a\}.
\]
In this work, all measures on Euclidean spaces will be assumed to be Borel probability measures which vanish on hyperplanes. 
Such measures are sometimes called \emph{mass distributions}.
\begin{definition}
	Let $N \geq 1$ be an integer. The \emph{depth of a point} $x$ with respect to the measure $\mu$ is
	\[
		\depth_{\mu}(x) := \inf \{ \mu( H^0_{v,a} ) \mid H_{v,a} \text{ is an oriented affine hyperplane with } x \in H^0_{v,a}  \}.
	\]
\end{definition}
In the literature there are some other notions of depth. So, in order to distinguish our notion from the others, we recall that the depth we are considering is also called half-space depth or Tukey depth \cite{tukey1975mathematics}. An important result concerning the depth of a point is the Rado theorem \cite{rado1946theorem}, which states that for every measure $\mu$ on $\mathbb{R}^d$ there exists a point $x$ such that $\depth_{\mu}(x) \geq \frac{1}{d+1}$. This result is also known as the \emph{centerpoint theorem}.

\medskip
Our interest in the study of the depth of a point, as well as its applications, comes from the following result obtained by 
Dol'nikov  in \cite{dol1987generalized, dol1992generalization} and independently by 
\v{Z}ivaljevi\'c \& Vre\'{c}ica in \cite{zivaljevic1990extension}: 
Let $\mu_1 , \ldots , \mu_m$ be $m$ measures in $\mathbb{R}^N$, where $m \leq N$. Then there is an $(m-1)$-dimensional affine subspace $L$ such that every half-space containing $L$ contains a fraction of at least $\frac{1}{N-m+2}$ of each measure.
This classical result, known as the \emph{center transversal theorem}, can be stated in terms of depth of a point as follows:

\begin{theorem}[Center transversal theorem]
\label{Thm CenterTT}
	Let $m$, $n$, and $N$ be positive integers with $N \geq m + n -1$. For every collection of $m$ measures 
	$\mu_1 , \ldots , \mu_m$ on $\mathbb{R}^N$, there exists an $n$-dimensional linear subspace $\Gamma$ and a point 
	$x \in \Gamma$ such that for every $1 \leq i \leq m$
	\[
		\depth_{\Gamma_* \mu_i}(x) \geq \frac{1}{n+1}.
	\]
\end{theorem}
Here $\Gamma_* \mu$ refers to the marginal measure of $\mu$ with respect to the subspace $\Gamma$ defined as follows:
\[
	\Gamma_* \mu(X) := \mu( \pi_{\Gamma}^{-1}(X) ),
\]
for every $X \subseteq \Gamma$, where $\pi_{\Gamma} \colon \mathbb{R}^N \to \Gamma$ denotes the orthogonal projection onto \(\Gamma\). Notice that $\pi_{\Gamma}^{-1}(x)$ is an $(N-n)$-dimensional affine subspace such that every half-space containing it contains a fraction of at least $\frac{1}{N-m+2}$ of each measure because $N \geq m + n -1$. 
Moreover, Theorem \ref{Thm CenterTT} also generalizes the Rado theorem, which is the case $m=1$.

\medskip
Recently, a version of the center transversal theorem with an improved bound on the depth was proved in \cite[Theorem 1.6]{blagojevic2018center}. Applying fairly advanced techniques of algebraic topology, the authors obtained the following result:

\begin{theorem}[Center transversal theorem with an improved Rado depth]
\label{Thm CTT}
	Let $m \geq 1$ and $k \geq 2$ be integers, and let 
		\begin{itemize}
			\item $d \geq 2m + k -1$ if $k +1$ is not a power of $2$, and
			\item $d \geq 3m + k -1$ if $k +1$ is a power of $2$.
		\end{itemize}
	For every collection of $m$ measures $\mu_1, \ldots , \mu_m$ on $\mathbb{R}^d$, there exists a $k$-dimensional linear 
	subspace $\Gamma$ and a point $x \in \Gamma$ with the property that all marginal measures 
	$\Gamma_* \mu_1, \ldots , \Gamma_* \mu_m$ satisfy
	\[
		\depth_{\Gamma_* \mu_i}(x) \geq \frac{1}{k+1} + \frac{1}{3( k+1 )^3}.
	\]
\end{theorem}

The term ``Rado depth'' comes from the bound on the depth obtained by Rado in \cite{rado1946theorem}, which, in particular, applies to every marginal measure as a lower bound, and the number in the Rado theorem is usually called the \emph{Rado bound}.

\medskip
Theorem \ref{Thm CTT} generalizes previous work of Magazinov \& P\'or \cite[Theorem 1]{magazinov2018improvement} and \cite[Theorem 1.4]{blagojevic2018center}, called the centerline theorem, which determines the depth of a point over one marginal measure. 
Theorem \ref{Thm CTT} also represents an important extension of the classical center transversal theorem in which the required depth is improved at the cost of increasing the dimension of the ambient space. This increase is linear in the number of measures and the dimension of $\Gamma$. In this paper, we will present an extension of Theorem \ref{Thm CTT} to mass assignments.
%

\section{ \Large An Improved Rado Depth for Mass Assignments}
\label{sec: Sec1}

We start by introducing mass assignments, which will take the place of the marginals of measures in Theorem \ref{Thm CTT}.

\subsection{Mass Assignments}
\label{subsec: MA}

Let $\mathcal{M}_+(X)$ be the space of all finite Borel measures on a topological space $X$ equipped with the weak topology. 
That is, the minimal topology such that for every bounded and upper semi-continuous function $f \colon X \to \mathbb{R}$, the induced function $\mathcal{M}_+(X) \to \mathbb{R}$, $\mu \mapsto \int f d\mu$, is upper semi-continuous. 
For a definition of a Borel measure on a topological space, see \cite[Def. 2.15]{Rudin-1987}. In case $X$ is an Euclidean space, the space of all mass distributions, or simply measures under our standing assumptions, will be denoted by $\mathcal{M}(X) \subseteq \mathcal{M}_+(X)$.

\medskip
Let $E$ be a real vector bundle over a path-connected space $B$ with fiber $E_b$ at $b \in B$. Let us consider the associated fiber bundle 
\begin{equation}\label{ma}
	\mathcal{M}(E) := \{ (b,\mu) \mid b \in B \text{ and } \mu \in \mathcal{M}(E_b)  \} \longrightarrow B
\end{equation}
given by $(b, \mu) \mapsto b$.
Any cross-section $\mu \colon B \to \mathcal{M}(E)$ of the fiber bundle \ref{ma} is called a \emph{mass assignment} on the Euclidean vector bundle $E$. Notice that $\mu^b := \mu(b)$ is a measure on the fiber $E_{b}$ for every $b \in B$. For a more detailed treatment of mass assignments on Euclidean vector bundles, see \cite{blagojevic2023many}.
The name ``mass assignment'' was chosen to align with ``mass distribution'' for the measures on the fibers; however, as stated in the introduction and following the terminology in \cite{blagojevic2018center}, we will use the term ``measure'' for ``mass distribution''. 

\medskip
Motivated by the study of the center transversal theorem for mass assignments in \cite{schnider2020ham} and \cite{axelrod2022bisections},  by the results obtained in \cite{blagojevic2018center} regarding the improved Rado depth of the 
measures,  and by the study of mass assignments over arbitrary Euclidean vector bundles in \cite{blagojevic2023many},  one can formulate the following general problem:

\emph{Given a Euclidean vector bundle $\eta \colon \mathbb{R}^k \hookrightarrow E \to B$, determine conditions under which, for every collection of $m$ mass assignments $\mu_1, \ldots , \mu_m$ on $\eta$, there exist points $b \in B$ and $x \in E_b \cong \mathbb{R}^k$ such that all the measures $\mu^{b}_1, \ldots , \mu^{b}_m$ on $E_b$ have sufficient depth with respect to the point $x$.}

(Here, ``sufficient depth'' means depth bounded below by some constant depending only on $B$.) This kind of problem has been studied using advanced techniques of algebraic topology, more precisely, using results concerning characteristic classes. Those techniques so far seem to require specific computations of characteristic classes of the bundle $\eta$, which means they can only be applied to a specific choice of vector bundle. Mass assignments on tautological vector bundles over Grassmannian manifolds have been recently used in \cite{schnider2020ham}, \cite{axelrod2022bisections}, and \cite{BlagojevicCalles-2023}. Those papers use this particular kind of mass assignment to study extensions of mass partition problems like the Gr\"unbaum--Hadwiger--Ramos problem and the center transversal theorem. We have chosen instead to focus on a slightly more general class of bundles that includes tautological bundles on Grassmannians. The bundles we consider are defined on real flag manifolds, and are described in the next section.

\subsection{Mass Assignments Over Flag Manifolds}
\label{subsec: MA Flag}

Let $n_1, \ldots , n_r$ be positive integers and let $n = n_1 + n_2 + \cdots + n_r$. The \emph{real flag manifold} 
$\Fl(n_1 , \ldots , n_r)$ is the set of $r$-tuples $(V_1, \ldots, V_r)$ of  vector subspaces of $\IR^n$ such that
\[
	V_1 \subset V_2 \subset \cdots \subset V_r = \IR^{n},
\]
with $\dim(V_i) = \sum_{j = 1}^{i} n_j$. Each element in $\Fl(n_1 , \ldots , n_r)$ will be represented by an $(r-1)$-tuple 
$(V_1, \ldots , V_{r-1})$, leaving out the subspace $V_r$ since it always refers to $\IR^{n}$.

\medskip
There is an equivalent alternative description of points on a real flag manifold, which is convenient for describing the cohomology ring of $\Fl(n_1 , \ldots , n_r)$; namely, we can describe a point on a real flag manifold as an $r$-tuple
$(W_1, \ldots , W_r)$ of subspaces of $\mathbb{R}^n$ which are mutually orthogonal and satisfy $\dim(W_i) = n_i$ and \(\bigoplus_{i=1}^r W_i = \IR^n\). It is easy to go back and forth between both representations: $V_i = W_1 \oplus \cdots \oplus W_i$, and \(W_i\) is the orthogonal complement of \(V_{i-1}\) inside \(V_i\).
This new description allows us to define a vector bundle $\omega_i$ of rank $n_i$ associated to $W_i$, namely:
\[
	\mathbb{R}^{n_i} \hookrightarrow E(\omega_i) \xrightarrow{\pi} \Fl(n_1 , \ldots , n_r) 
\]
where
\[
	E(\omega_i) = 
		\big\{ (W_1,W_2, \ldots , W_r , v) \in \Fl(n_1 , \ldots , n_r) \times \IR^{n} \mid v \in  W_i  \big\},
\]
and  $\pi(W_1,W_2, \ldots , W_r , v) = (W_1,W_2, \ldots , W_r)$. The vector bundle $\omega_i$ is called the 
\emph{$i$-th tautological vector bundle} over $\Fl(n_1 , \ldots , n_r)$. Notice that $\omega_1 \oplus \cdots \oplus \omega_n$ is a trivial vector bundle. 
Also, since $\Fl(k , n-k) = \Gr{k}{n}$, the associated bundle $\omega_1$ is actually the tautological vector bundle over the Grassmannian $\Gr{k}{n}$, usually denoted by $\gamma_{k}^{n}$.
For more details about flag manifolds, see \cite[Section 9.5]{hausmann2014mod}.

\medskip
The following result describes the cohomology ring of $\Fl(n_1 , \ldots , n_r)$ in terms of the Stiefel-Whitney classes of the tautological bundles $\omega_i$ introduced above. 
\begin{theorem}\label{thm cohom Flag}
	The cohomology ring $H^*( \Fl(n_1,n_2, \ldots , n_r) ; \IF_2 )$ is isomorphic to the quotient of the polynomial ring
	\[ 
		 \IF_2[w_i( \omega_j )], \hspace{.5cm} 1 \leq i \leq n_j \hspace{.5cm} \text{with} \hspace{.5cm} j = 1, \ldots , r
	\]
	by the ideal generated by the homogeneous components of the total Stiefel-Whitney classes 
	$w(\omega_1) \cdots w(\omega_r)$ in positive degrees.	
\end{theorem}
For details about the proof, see \cite[Theorem 9.5.14]{hausmann2014mod}. From now on, we will be using the nested subspaces interpretation of a real flag manifold. Notice that the result obtained in \ref{thm cohom Flag} for $\Fl(k , n-k) = \Gr{k}{n}$ coincides with the classical result of Borel \cite{borel1953cohomologie}.

\medskip
Having defined tautological vector bundles over flag manifolds,  we will keep our attention on mass assignments regarding such bundles.  More precisely,  we will be referring to \emph{mass assignments on the vector bundle 
$\nu_i :=  \bigoplus_{j=1}^{i} \omega_j$ over the flag manifold $\Fl(n_1 , \ldots , n_r)$} as a cross-section of the bundle $\nu_i$.
Due to the equivalence between the definitions of flag manifold, we will also be referring to $\nu_i$ as a tautological vector bundle,  for $1 \leq i \leq r-1$. We will always say explicitly which bundle we mean, so there is no risk of confusion from calling all these bundles tautological.


\medskip
We shall focus on the following problem:

\begin{problem}
\label{PrincipalProblem}
	Determine all tuples of positive integers $(m,  i,  n_1,  n_2, \ldots , n_r)$ such that for every collection of $m$ 
	mass assignments $\mu_1, \ldots , \mu_m$ on  $\nu_i$ over $\Fl(n_1 , \ldots , n_r)$,  with $ 1 \leq i \leq r-1$,  there exists a flag 
	\[
			\mathcal{F} = ( V_{1},  V_{2},  \ldots , V_{{r-1}} ) \in  \Fl(n_1, n_2, \ldots , n_r),
	\]
	with $V_1 \subset V_2 \subset \cdots \subset V_r = \IR^{n}$ and $\dim(V_s) = \sum_{j = 1}^{s} n_j$, and a point 
	$x \in V_i$ with the property that all the measures $\mu^{\mathcal{F}}_1, \ldots , \mu^{\mathcal{F}}_m$ 
	on $V_i$ have sufficient depth with respect to the point $x$.
\end{problem}

\medskip 
Our strategy to study Problem \ref{PrincipalProblem} is to adapt the techniques used in \cite{blagojevic2018center} to the case of mass assignments over flag manifolds.  This implies replacing the tautological bundle over the Grassmannian $\Gr{k}{\ell}$ used in \cite{blagojevic2018center} by the vector bundle $\nu_i$ over the real flag manifold $\Fl(n_1, n_2, \ldots , n_r)$.

\subsection{The Center Transversal Theorem for Mass Assignments}
\label{subsec: General results}

We will present a version of \cite[Theorem 2.4]{blagojevic2018center} using flag manifolds. Such result is a partial answer for Problem \ref{PrincipalProblem} and will be proved in Section \ref{sec: Proof}. To be more precise, we shall prove the following result:

\begin{theorem}
\label{Principal thm}
	Let $m, n_1, \ldots , n_r$ be positive integers, with $r \geq 3$ and $n_{i+2} + \cdots + n_r \geq 2m-1$, and let $1 \leq i \leq r-2$. 
	For every collection of $m$ mass assignments $\mu_1, \ldots , \mu_m$ on $\nu_i$ over $\Fl(n_1 , \ldots , n_r)$, there 
	exists a flag $ \mathcal{F} = ( V_{1}, V_{2}, \ldots , V_{{r-1}} ) \in \Fl(n_1, n_2, \ldots , n_r) $, and a point $x \in V_{i}$ such that for 
	every $1 \leq j \leq m$:	
	\[
		\depth_{\mu^{\mathcal{F}}_j }(x) \geq \frac{1}{k+1} + \frac{1}{3( k+1 )^3},
	\]
	with $k = \dim(V_i)$.
\end{theorem}

There is a special case of Theorem \ref{Principal thm} that we think is important to highlight: 
Let us start by considering a collection of mass assignments on $\omega_1$ over $\Fl(k,\ell - k, d - \ell)$, where $k$ and $\ell$ are positive integers, with $k < \ell$ and $d - \ell \geq 2m -1$. Moreover, the mass assignments under consideration will be marginals of mass assignments on $\Gr{\ell}{d}$, that is, 
\[
	\mu^{(\Gamma, L)}_i := \Gamma_* \mu^{L}_i.
\]
All these considerations lead us to the following result:

\begin{corollary} 
\label{Corol 1} 
	Let $m$, $k$, and $\ell$ be positive integers with $k < \ell$, and let $d \geq 2m + \ell -1$.
	For every collection of $m$ mass assignments $\mu_1, \ldots , \mu_m$ on $\gamma_{k}^{d}$ over $\Gr{k}{d}$, 
	there exists a $k$-dimensional linear subspace $\Gamma$ contained in $L \in \Gr{k}{d}$, and a point $x \in \Gamma$ 
	such that for every $1 \leq i \leq m$:
	\[
		\depth_{\Gamma_*\mu^{L}_i}(x) \geq \frac{1}{k+1} + \frac{1}{3( k+1 )^3}.
	\]
\end{corollary}

Corollary \ref{Corol 1} is a variant of \emph{the center transversal theorem with an improved Rado depth} \cite[Theorem 1.6]{blagojevic2018center} over linear subspaces $L \in \Gr{\ell}{d}$, where  $\ell$ is lower than the dimension obtained in \cite{blagojevic2018center}. 
Notice also that, in the same way as the improvement of the Rado depth presented in \cite{blagojevic2018center} required a larger dimension of the corresponding Euclidean space, our improvement on the dimension $\ell$ is reflected in a larger dimension for the ambient space in the Grassmannian manifold. 

\medskip
Our theorem does not include the case $k=\ell$, which would have been a generalization of \cite[Theorem 1.6]{blagojevic2018center} to mass assignments over a Grassmannian, and without a case distinction based on whether $\ell+1$ is a power of $2$ or not. We do not see how to avoid the case distinction in that case, but the following generalization of \cite[Theorem 1.6]{blagojevic2018center} to mass assignments does hold simply by replacing the marginals $\Gamma_* \mu_i$ by $\mu_i^{\Gamma}$ throughout the proof in \cite[Section 4]{blagojevic2018center}.

\begin{proposition}[Center Transversal theorem for mass assignments with an improved Rado depth] 
\label{Prop1} 
	Let $m \geq 1$, $k \geq 2$ be integers, and let 
		\begin{itemize}
			\item $d \geq 2m + k -1$ if $k +1$ is not a power of $2$, and
			\item $d \geq 3m + k -1$ if $k +1$ is a power of $2$.
		\end{itemize}
	For every collection of $m$ mass assignments $\mu_1, \ldots , \mu_m$ on $\gamma_{k}^{d}$ over $\Gr{k}{d}$,  there exists an 
	$k$-dimensional linear subspace $\Gamma$ and a point $x \in \Gamma$ such that for every $1 \leq i \leq m$:
	\[
		\depth_{ \mu^{\Gamma}_i }(x) \geq \frac{1}{k+1} + \frac{1}{3( k+1 )^3}.
	\]
\end{proposition}

\section{Proof of theorem \ref{Principal thm}}
\label{sec: Proof}

The proof of theorem \ref{Principal thm} is based on the one given by \cite{blagojevic2018center} in section $4$.
Before presenting the proof of our main result we need to mention some considerations.
We will be working with the \emph{depth of a measure} $\mu$ on $\mathbb{R}^n$ defined as follows:
\[
	\depth(\mu) := \text{sup}_{x \in \mathbb{R}^n} \depth_{\mu}(x).
\]
To choose the point $x \in V_i$ we will rely on a continuous function constructed in \cite{blagojevic2018center} which assigns to each measure $\mu$ on $\mathbb{R}^n$ a point $c(\mu) \in \mathbb{R}^n$ where,  roughly speaking, the depth of the measure is maximized. More precisely, \(c(\mu)\) satisfies:
\[
	\depth(\mu) < \frac{1}{n+1} + \frac{1}{3( n+1 )^3} \hspace{.5cm} \Longleftrightarrow \hspace{.5cm} 
		\depth_{\mu}( c(\mu) ) < \frac{1}{n+1} + \frac{1}{3( n+1 )^3}.
\]
This is all we will need about \(c(\mu)\), but for a more detailed explanation of the construction and properties of the point $c(\mu)$ see \cite{magazinov2018improvement} and \cite{blagojevic2018center}.

\medskip
We seek to prove that for every collection of $m$ mass assignments $\mu_1, \ldots , \mu_m$ on $\nu_i$ over 
$\Fl(n_1,  n_2 , \ldots , n_r)$, there exists a flag $\mathcal{F}=(V_1, V_2, \ldots , V_{r-1})$ in $\Fl(n_1,  n_2 , \ldots , n_r)$ such that for every $1 \leq j \leq m$	
\[
	\depth( \mu^{\mathcal{F}}_j  ) \geq \frac{1}{k+1} + \frac{1}{3( k+1 )^3} \; ,
\]
and in addition
\[
	c( \mu^{\mathcal{F}}_1 ) = \ldots = c( \mu^{\mathcal{F}}_m ).
\]
Moreover,  since all masses $\mu^{\mathcal{F}}_j$ are defined on $V_i$,  our choice of $x$ will be the point 
$c( \mu^{\mathcal{F}}_1 ) = \ldots = c( \mu^{\mathcal{F}}_m )$ in $V_i$.

\medskip 
Analogously to \cite[Section 4]{blagojevic2018center}, we consider the following open sets in the real flag manifold:
\[
	U_0 = \Big\{ \mathcal{F} \in \Fl(n_1,  n_2 , \ldots , n_r) \mid 
		c( \mu^{ \mathcal{F} }_{1} ), \ldots , c( \mu^{ \mathcal{F} }_{m} ) \text{ do not all coincide} \Big\}
\]
and 
\[
	U_j = \Big\{ \mathcal{F}  \in \Fl(n_1,  n_2 , \ldots , n_r) \mid 
		\depth( \mu^{ \mathcal{F} }_j ) < \frac{1}{k+1} + \frac{1}{3(k+1)^3}  \Big\}
\]
for $1 \leq j \leq m$. Any point outside of \(U_0 \cup U_1 \cup \cdots \cup U_m\) satisfies the conclusion of Theorem \ref{Principal thm}, so we shall assume that \(\{U_j : 0 \le j \le m\}\) is an open cover of the flag manifold and obtain a contradiction, To obtain this contradiction, still following the strategy of \cite[Section 4]{blagojevic2018center}, we will define sections of certain vector bundles over each open set of the cover and then use \cite[Lemma 3.2]{blagojevic2018center}, which is reproduced below, coincidentally also as Lemma \ref{lemma: cohom class}, and whose proof we include for convenience of the reader. The proof relies on the following standard result, whose proof we include for the same reason.

\begin{lemma}
\label{lemma mult classes}
	Let $\{ U_i \}_{i\in I}$ be a cover of $X$, and let $\alpha_1, \ldots , \alpha_n$ be cohomology classes in $H^*(X; \IF_2)$. 
	Consider the inclusion maps $j_i \colon U_i \to X$. If ${j_i}^*(\alpha_i)=0$ for every $i \in I$, then $\alpha_1 \cdots \alpha_n=0$.
\end{lemma}
\begin{proof}
	Since ${j_i}^*(\alpha_i)=0$ for every $i \in I$,  by the long exact sequence 
	\[
		\cdots \rightarrow \tilde{H}^{n}( X/U_i ;\IF_2) \xrightarrow{{q_i}^*} \tilde{H}^{n}( X ;\IF_2) \xrightarrow{{j_i}^*}
		\tilde{H}^{n}( U_i ;\IF_2) \rightarrow \cdots
	\]
	there is a class $\beta_i \in \tilde{H}^{n}( X/U_i ;\IF_2)$ such that ${q_i}^*(\beta_i) = \alpha_i$.  Consider now the commutative 
	diagram 	
	\[
		\xymatrixcolsep{.8in}
		\xymatrix{
			X \ar[r]^{\Delta} \ar[rd]_{\mathrm{const}} & X \wedge \cdots \wedge X \ar[d]^{q_1 \wedge \cdots \wedge q_n} \\
			                         & X/U_1 \wedge \cdots \wedge X/U_n
	}\]	
	where $\Delta$ is the diagonal map.  Finally,  since 
	\[
		\Delta^* \circ (q_1 \wedge \cdots \wedge q_n)^* ( \beta_1 \wedge \cdots \wedge \beta_n ) = \alpha_1 \cdots \alpha_n,	
	\]
	and on the other hand $\Delta^* \circ (q_1 \wedge \cdots \wedge q_n)^*$ is identically zero,  then $\alpha_1 \cdots \alpha_n	 = 0$.
\end{proof}

	\begin{lemma}{\cite[Lemma 3.2]{blagojevic2018center}}
		\label{lemma: cohom class}
		Let $\pi : E \to X$ be a fiber bundle, $f : B \to X$ any map, and $B = \bigcup_{i=1}^m U_i$ be an open cover. Assume there is a cohomology class $\alpha \in H^*(X; \IF_2)$ such that $\pi^*(\alpha) = 0$ and ${f^*(\alpha)}^m \neq 0$. Then there cannot exist continuous local sections of the pullback bundle $f^*\pi$ over each of the open sets $U_1, \ldots , U_m$.
	\end{lemma}
	\begin{proof}
          The proof of the lemma is by contradiction. Let us suppose that there is a section $s_i : U_i \to E(f^{*}\pi)$ of $f^*\pi$ for each $i=1, \ldots, m$. Since $(f^{*}(\alpha))^m \neq 0$, by Lemma \ref{lemma mult classes}, it must be the case that for some $j$ we have ${f|_{U_j}}^* (\alpha) \neq 0$. Consider now the following commutative square:
		\[
			\xymatrix{
							 E(f^*\pi)|_{U_j}   \ar[r]   \ar[d]^{p}  & E(f^*\pi)   \ar[r]^{\hat{f}}   \ar[d]   &   E  					
							 \ar[d]^{\pi}
							 \\
							 U_j  \ar@{^(->}[r] \ar@/^/[u]^{s_j} & B   \ar[r]^{ f } & X.
			}		
		\] 
		Since $p \circ s_j = \id_{U_j}$, we have $s_j^{*} \circ p^{*} = \id$, and thus $p^*$ is injective.  Consequently, 
		\[  {\hat{f}|_{U_j}}^* \circ \pi^* (\alpha)  = p^* \circ {f|_{U_j}}^* (\alpha) \neq 0,\] contradicting the assumption that $\pi^{*}(\alpha) = 0$.
	\end{proof}

        With the general results powering the strategy out of the way, we now turn back to specifics. Our first task is to describe the bundle and sections to which we shall apply Lemma~\ref{lemma: cohom class}. For that purpose, we will present a description of the Stiefel manifold $ V_{k}(\mathbb{R}^{\infty})$ in terms of regular unit simplices with ordered vertices and whose center is the origin. The space of such simplices contained in fibers of the tautological bundle $\gamma_k^{\infty}$ can be parametrized by the Stiefel manifold $ V_{k}(\mathbb{R}^{\infty})$ as follows: first, identify $ V_{k}(\mathbb{R}^{\infty})$ with the space of all orthogonal linear maps $\mathbb{R}^k \to \mathbb{R}^{\infty}$; then, fix a model unit centered regular $k$-simplex $\Delta_k \subset \mathbb{R}^k$ and make the correspondence $\phi \mapsto \phi(\Delta_k)$, where $\phi \in V_{k}(\mathbb{R}^{\infty})$.

There is a fiberwise free action of the symmetric group $\mathcal{S}_{k + 1}$ on the Stiefel manifold $ V_{k}(\mathbb{R}^{\infty})$ by permuting the vertices of the regular simplices. In this way, $ V_{k}(\mathbb{R}^{\infty}) / \mathcal{S}_{k + 1}$ parametrizes the space of all centered regular unit simplices with unordered vertices contained in the fibers of the tautological bundle $\gamma_k^{\infty}$. Moreover, since the infinite Stiefel manifold is contractible, $ V_{k}(\mathbb{R}^{\infty}) / \mathcal{S}_{k + 1}$ is a model for $B\mathcal{S}_{k + 1}$.

\medskip
In addition to this,  $B\mathcal{S}_{k+1}$ can be seen as the total space of the following fiber bundle 
\[
		\eta \colon  O(k)/\mathcal{S}_{k+1} \to B\mathcal{S}_{k+1} \xrightarrow{\sigma} BO(k)
\] 
where the projection map $\sigma$ is induced by the representation $ \rho \colon \mathcal{S}_{k+1} \to O(k)$ obtained from the permutations of the vertices of the regular simplex. Consider now the following pullback square and its associated diagram in cohomology

\begin{center}
  \begin{tabular}{ c c c  }
          $
                  \xymatrix{
                           E(f^*\eta)   \ar[r]   \ar[d]^{\sigma}    &   B\mathcal{S}_{k+ 1}  \ar[d]^{\sigma_k} \\
                           B   \ar[r]^{ f } & BO(k)
                  }	 	
          $
          & \hspace{1cm} &
          $
                  \xymatrix{
                           H^*(E(f^*\eta; \IF_2)    &   H^*(B\mathcal{S}_{k+ 1}; \IF_2)  \ar[l]  \\
                           H^*(B; \IF_2) \ar[u]_{\sigma^*}    & H^*(BO(k); \IF_2) \ar[l]_{ f^* }  \ar[u]_{\sigma_k^*}.
                  }	 	
          $
  \end{tabular}
\end{center}
              
What follows is to define local sections of this bundle over the open cover $\{U_1, \ldots , U_m\}$ of the complement of $U_0$ in the flag manifold. Notice that a section represents a continuous choice of a centered unit regular simplex contained in each fiber of the bundle.
		
Let $1 \leq j \leq m$ and a consider a point $\mathcal{F} := (V_1, \ldots, V_{r-1}) \in U_j$. Then, $\mu^{ \mathcal{F} }_j$ has an associated regular simplex $\Delta( \mu^{ \mathcal{F} }_j) \subset V_i$ centered at the origin, which is described in \cite[Section 2]{blagojevic2018center}, so we can define a map on $U_i$ via the assignment
	\[
		\mathcal{F} \longmapsto (\Delta(\mu^{ \mathcal{F} }), V_i^\perp),
	\]
where $V_i^\perp$ denotes the orthogonal complement of $V_i$ inside $V_{i+1}$. This map is a section of the following pullback bundle:
\begin{equation} \label{diagram: pullback}
  \resizebox{.76\textwidth}{!}{
    \xymatrix{
      E   \ar[r]   \ar[d]    &   E'  \ar[r]   \ar[d]      &      E( \rho^*(\sigma) )   \ar[r]   \ar[d]  & 
      B\mathcal{S}_{k+ 1} \times BO(\ell -k) 
      \ar[d]^{\sigma} \\
      U_i \ar[r] \ar@/^/[u] & \Fl(n_1,  n_2 , \ldots , n_r) \ar[r]^{\; \; \pi} & \Gr{\ell}{n} \ar[r]^{\rho} & BO(\ell),
    }
  }
\end{equation}
where
\[
 k = \text{dim}(V_i) = \sum_{s=1}^{i} n_s \hspace{.5cm} \text{and} \hspace{.5cm} \ell = \text{dim}(V_{i+1}) = k + n_{i+1},
\]
the projection map $\sigma$ is the composite
	\[
		B\mathcal{S}_{k+1} \times BO(l-k) \rightarrow BO(k) \times BO(\ell-k) \rightarrow BO(\ell),
	\]
and the map $B\mathcal{S}_{k+1} \to BO(k)$ is induced by the inclusion of the symmetries of a regular simplex.

We now have a bundle defined over the complement of $U_0$ and local sections of that bundle defined on $U_1, U_2, \ldots, U_m$, which we have assumed form an open cover of $U_0^c$ in order to obtain contradiction by applying Lemma~\ref{lemma: cohom class}. The final ingredient needed to apply the lemma is a cohomology class on the base of the bundle.  We will see now that the Stiefel-Whitney class $w_\ell$ in $H^*(BO(\ell); \IF_2)$ satisfies the hypothesis of the lemma. We need to show that $\sigma^{*}(w_{\ell}) =0$ and $f^{*}(w_{\ell})^m\neq0$, where $f = (\rho \circ \pi)|_{U_0^c}$. The first condition is straightforward, because $\sigma$ factors through $BO(k) \times BO(\ell-k)$ and we have $k<\ell$ and $\ell-k<\ell$.

\medskip

It only remains to show that $f^{*}(w_{\ell})^m\neq0$. Consider now the class $w_{\ell}^{n-\ell}$ in $H^*( \Gr{\ell}{n} ; \IF_2 )$,  which by \cite[Lemma 1.2]{hiller1980cohomology}, is not trivial.  Since the projection map 
\[
	\pi \colon \Fl(n_1,  n_2 , \ldots , n_r)  \longrightarrow  \Gr{\ell}{n}
\]
induces a monomorphism in cohomology,  and $\pi^*(w_{\ell}) = w_{n_1}(\omega_1) \cdots w_{n_{i+1}}(\omega_{i+1})$,  we have
\[
	\pi^*(w_{\ell}^{n-\ell}  ) = \big(  w_{n_1}(\omega_1) \cdots w_{n_{i+1}}(\omega_{i+1})  \big)^{n-\ell} \neq 0 \quad (\text{on } \Fl(n_1, n_2, \ldots, n_r)).
\]  

Next we claim that ${ \big( w_{n_1}(\omega_1) \cdots w_{n_i}(\omega_i) \big) }^{m-1} = 0$ in $H^*(U_0; \IF_2)$. To see this, recall that for $\mathcal{F} \in U_0 $ the points $c( \mu^{ \mathcal{F} }_{1} ), \ldots , c( \mu^{ \mathcal{F} }_{m} )$ do not all coincide, which implies that the function
\[
  \mathcal{F} \longmapsto \Big( c( \mu^{ \mathcal{F} }_{j+1} ) - c( \mu^{ \mathcal{F} }_{j} ) \Big)_{j=1}^{m-1}
\]
defines a nonzero section of the restriction of the Whitney sum $\nu_i ^{\oplus m-1} $ to $U_0$. Now, the top Stiefel-Whitney class represents an obstruction to the existence of nonzero sections on a vector bundle \cite[Property 9.7]{milnor1974characteristic}, and thus
\[
  {\big( w_{n_1}(\omega_1) \cdots w_{n_i}(\omega_i) \big) }^{m-1} = 0 \quad (\text{on } U_0).
\]
If also we had
\[
  \big(  w_{n_1}(\omega_1) \cdots w_{n_{i+1}}(\omega_{i+1})  \big)^{n-\ell-m+1} = 0 \quad(\text{on } U_0^c),
\]
then, by Lemma~\ref{lemma mult classes}, we would have $\pi^{*}(w_\ell^{n-\ell})=0$, since $\pi^{*}(w_\ell^{n-\ell})$ is a multiple of the product of the classes we know are zero when restricted to $U_0$ and $U_0^c$. But we have seen that $\pi^{*}(w_\ell^{n-\ell})\neq0$, and therefore $f^*(w_{\ell})^{n-\ell-m+1} \neq 0$.  By assumption,  
$m \le n - \ell - m + 1$, so we also must have $f^{*}(w_\ell)^m\neq 0$ as required.

\subsection*{Acknowledgements}

The second author was supported by the UNAM Posdoctoral Scholarship Program (DGAPA).


\bibliography{Ref.bib}{}
\bibliographystyle{amsplain}

\end{document}